\documentclass[10pt, reqno, a4wide]{amsart}
\usepackage{a4wide}
\usepackage{amsfonts,amsthm}
\usepackage{graphicx} 
\usepackage{xcolor}
\usepackage{url}
\usepackage{enumitem}
\usepackage{hyperref}
\usepackage[left=1in,right=1in,top=0.6in,bottom=1in]{geometry}
\DeclareRobustCommand{\van}[3]{#2} 

\theoremstyle{plain}

\title{An elementary proof of the bunkbed conjecture for forests}
\author{Serte Donderwinkel$^1$, Joost Jorritsma$^2$, Guillem Perarnau$^{3,4}$}
\address{$^1$Bernoulli Institute of Mathematics, University of Groningen, $^2$Department of Statistics, University of Oxford, $^3$Universitat Polit\`ecnica de Catalunya (UPC), Barcelona, Spain,
$^4$Centre de Recerca Matem\`atica, Bellaterra, Spain.
}
\email{s.a.donderwinkel@rug.nl, joost.jorritsma@stats.ox.ac.uk, guillem.perarnau@upc.edu}
\newcommand{\Prob}{\mathbb{P}}

\begin{document}
\begin{abstract}
   Although false for general graphs, this note gives an elementary  proof of the bunkbed conjecture for any acyclic graph. The argument is short and self-contained, and may be of educational interest.
\end{abstract}
\maketitle

\vspace{-0.5cm}
Let $G^+, G^-$ be two copies of a finite graph $G$ with vertex labels $[n]:=\{1,\ldots, n\}$. We consider these graphs as stacked on top of each other, and let the vertices be labeled $[n]^+:=\{1^+,\ldots, n^+\}$ and $[n]^-:=\{1^-,\ldots, n^-\}$. 
Fix $H\subseteq [n]$, and let $G^\pm=(V^\pm,E^\pm)$ be the graph with $V^\pm=V^+\cup V^-$ and 
\[
E^\pm := E^+\cup E^-\cup \big\{\{v^+, v^-\}: v\in H\big\}.
\]
Let $\Prob^G_{H,p}$ denote the law of the graph $G^\pm$ in which each edge in $E^+\cup E^-$ is open independently with probability $p\in(0,1)$, and all edges of the form $\{v^+, v^-\}$ are open.  We write $\{u\leftrightarrow v\}$ if there exists an open path from $u$ to $v$ in this graph. Kasteleyn conjectured (see~\cite{vdBK2001})

\medskip
\noindent\textbf{Former bunkbed conjecture.}
\emph{
    For any graph $G=(V,E)$, $p\in(0,1)$, and  $u,v\in [n]$, 
    \begin{equation}\label{bunkbed}
    \Prob_{H,p}^G\big(u^+\leftrightarrow v^+\big)\ge  \Prob_{H,p}^G\big(u^+\leftrightarrow v^-\big).
    \end{equation}
}

\medskip 
Intuitively, it should be more likely to maintain an open path between vertices on the same level, than to cross layers.
However, a counterexample was given recently by Gladkov, Pak and Zimin \cite{GPZ2025}, based on a counterexample of Hollom for the hypergraph version of the conjecture~\cite{H2025}. Since the conjecture fails in full generality, it becomes even more interesting to understand for which cases it remains valid. The conjecture has proven to be true for certain classes of graphs such as complete graphs~\cite{vHL2019}, wheels~\cite{L2009} or highly symmetric graphs~\cite{R2022} and when $p$ tends to $0$ or to $1$~\cite{H2024,HKN2023}. In particular, Linusson proved the conjecture for outerplanar graphs \cite{L2011}, which includes our setting here. We present an elementary proof of

\medskip 
\noindent\textbf{Theorem.}
\emph{
    The bunkbed conjecture holds for any forest (acyclic graph).
}

\medskip\noindent
Our theorem holds for any set $H$, and thus the theorem holds also for cases in which $H$ is random, independent of the percolation configuration. 
Another proof of our result is implied by~\cite{MP2024} which studies graphs `glued' along a vertex. We note that the bunkbed conjecture is false for some \emph{directed} acyclic graphs~\cite{przybylowski2025acyclic}.

We first prove~\eqref{bunkbed} for $G=P_n$ a path of $n-1$ edges, with consecutive vertices connected by an edge, and only considering $u=1$ and $v=n$. We then extend it to arbitrary forests, $u$ and $v$, by reduction to paths.

\begin{proof}[Proof for the path $P_n$]
    If $H$ is empty, the inequality is trivial, so we assume $H\neq\emptyset$. Let $m=\max H$. If there is a path from $1^+$ to $n^+$, there must be a path from $1^+$ to $m^+$, after which all remaining $n-m$ edges must be open. Analogous reasoning holds for paths to $n^-$. Since $m^+$ has an open edge to $m^-$, 
    \[
    \Prob_{H,p}^{P_n}\big(1^+\leftrightarrow n^+\big)=\Prob_{H,p}^{P_n}\big(1^+\leftrightarrow m^+\big)p^{n-m}=\Prob_{H,p}^{P_n}\big(1^+\leftrightarrow m^-\big)p^{n-m}=\Prob_{H,p}^{P_n}\big(1^+\leftrightarrow n^-\big).\qedhere
    \]
\end{proof}

\begin{proof}[Proof for arbitrary forests]
    Let $F$ be a forest, $u,v\in V(F)$ and  $H\subseteq V(F)$. Without loss of generality, we assume that $u$ and $v$ are in the same connected component of $F$, and that the consecutive vertices on the unique path from $u$ to $v$ have labels $(u=1,2,\ldots, \ell=v)$. Let $\Sigma_{[\ell]^c}$ be the percolation configuration restricted to the edges outside $[\ell]^\pm$. We show that the inequality \eqref{bunkbed} holds conditional on the event $\Sigma_{[\ell]^c}=\sigma_{[\ell]^c}$ for any $\sigma_{[\ell]^c}$, so that the unconditioned inequality follows by the tower property. For any $\sigma_{[\ell]^c}$, we will construct  $H'\subseteq[\ell]$ such that 
    \begin{align*}
    &\Prob_{H,p}^F\big(1^+\leftrightarrow \ell^+\mid \Sigma_{[\ell]^c}=\sigma_{[\ell]^c}) = \Prob_{H',p}^{P_\ell}\big(1^+\leftrightarrow \ell^+\big),   \text{ and }\\
&\Prob_{H,p}^{F}\big(1^+\leftrightarrow \ell^-\mid \Sigma_{[\ell]^c}=\sigma_{[\ell]^c}\big) = \Prob_{H',p}^{P_\ell}\big(1^+\leftrightarrow \ell^-\big).
    \end{align*}
    We let $H'$ be the vertices $z\in [\ell]$ for which there exists a path from $z^+$ to $z^-$ that does not use edges on the paths $[\ell]^\pm$, only traversing open edges in $\sigma_{[\ell]^c}$ and edges $\{w^+,w^-\}$ with $w\in H$. Thus, $H\cap[\ell]\subseteq H'$. Now, for any  configuration $\sigma_{[\ell]}$ on the edges in $[\ell]^\pm$, the following are equivalent:
    \begin{enumerate}[leftmargin=*]
    \item[(a)]there exists an open path from $1^+$ to $\ell^\pm$ using edges open in $\sigma_{[\ell]}$ or $\sigma_{[\ell]^c}$ and edges $\{w^+,w^-\}$ with $w\in H$; 
    \item[(b)]there exists an open path from $1^+$ to $\ell^\pm$  using only edges open in $\sigma_{[\ell]}$ and edges $\{w^+,w^-\}$ with $w\in H'$. 
    \end{enumerate}
    Thus, the claimed equalities follow, so that the theorem follows from the special case $P_\ell$. 
\end{proof}

\newpage

\newgeometry{margin=0.9in}

\noindent\textbf{Acknowledgements.\ }This note was written while the authors were attending the workshop `Combinatorics and Discrete Probability' at CIRM, Marseille, November 2025. JJ thanks Magdalen College for a Senior Demyship. SD acknowledges the financial support of the CogniGron research center and the Ubbo Emmius Funds (University of Groningen). Her research was also partially supported by the Marie Skłodowska-Curie grant GraPhTra (Universality in phase transitions in random graphs), grant agreement ID 101211705. GP acknowledges PID2023-147202NB-I00 funded by MICIU/AEI/10.13039/501100011033, and MSCA-RISE-2020-101007705.


\end{document}